\documentclass[10pt,a4paper]{aocarusodmiunict}
\issueinfo{}{}{}{}
\pagespan{1}{}
\PII{\textsc{}}
\copyrightinfo{2002}{{\textsc{Author: Andrea O. Caruso -- Date: June 2007}}}

\usepackage{amsmath,amsfonts,amssymb,hyperref,color,graphicx}
\usepackage[english]{babel}
\usepackage{pdfsync} 

\newtheorem{theorem}{Theorem}
\newtheorem*{theorem*}{Theorem}
\newtheorem{lemma}{Lemma}
\newtheorem*{lemma*}{Lemma}
\newtheorem{proposition}{Proposition}
\newtheorem*{proposition*}{Proposiion}
\newtheorem{corollary}{Corollary}
\newtheorem*{corollary*}{Corollary}
\newtheorem{definition}{Definition}
\newtheorem*{definition*}{Definition}
\newtheorem{notation}{Notation}
\newtheorem*{notation*}{Notation}
\newtheorem{remark}{Remark}
\newtheorem*{remark*}{Remark}
\newtheorem{example}{Example}
\newtheorem*{example*}{Example}

\newtheorem*{note*}{Note}
\newtheorem*{fremark*}{Concluding Remark}

\newcommand{\B}{\;\;\;\Box}

\newcommand{\ro}{\rho}

\newcommand{\tetaxe}{\Theta(\xi,\eta)}
\newcommand{\tetaex}{\Theta(\eta,\xi)}
\newcommand{\erre}{\mathbb{R}}
\newcommand{\errek}{{\mathbb{R}^k}}
\newcommand{\erren}{{\mathbb{R}^n}}
\newcommand{\erreN}{{\mathbb{R}^N}}

\newcommand{\erreq}{{\mathbb{R}^q}}
\newcommand{\ep}{\epsilon}

\newcommand{\paj}{\partial_{j}}

\newcommand{\gigno}{\mathbb{G}}
\newcommand{\gia}{\mathfrak{g}_{\,q,\,r}}
\newcommand{\gig}{\mathbb{G}_{\,q,\,r}}
\newcommand{\aij}{a^{ij}}

\newcommand{\aijtilde}{\widetilde{a^{ij}}}
\newcommand{\ajitilde}{\widetilde{a^{ji}}}
\newcommand{\xconi}{X_i}
\newcommand{\xconj}{X_j}

\newcommand{\yj}{Y_j}

\newcommand{\xjt}{X_j^T}

\newcommand{\xjttilde}{\widetilde{X_j}^T}
\newcommand{\xjtilde}{\widetilde{X_j}}

\newcommand{\xitilde}{\widetilde{X_i}}
\newcommand{\yuyq}{Y_1,\ldots,Y_q}
\newcommand{\xuxq}{X_1,\ldots,X_q}
\newcommand{\xuxqtilde}{\widetilde{X_1},\ldots,\widetilde{X_q}}
\newcommand{\xtilde}{\widetilde{X}}
\newcommand{\utilde}{\widetilde{u}}
\newcommand{\fcaptilde}{\widetilde{F}}
\newcommand{\fcapjtilde}{\widetilde{F^{j}}}
\newcommand{\lux}{\xjt(\aij\xconi)}
\newcommand{\luxtilde}{\xjttilde(\aijtilde\xitilde)}

\begin{document}

   \title[$W^{1,p}_{X}$ interior regularity
   with discontinuous coefficients]
   {$W^{1,p}_{X}$ interior estimates for\\
   variational hypoelliptic operator\\
   with $VMO_X$ coefficients}

   \author[A.O.\,Caruso]{A.O.\,Caruso}
   \address[]{Dipartimento di Matematica e Informatica
   \\ Universit\`a di Catania
   \\ Viale A. Doria 6--I, 95125, Catania, Italy}
   \email{\href{mailto:aocaruso@dmi.unict.it}{aocaruso@dmi.unict.it}}




   \keywords{Hypoelliptic operators, $L^{p}$ regularity,
   VMO spaces, Carnot--Caratheodory metric.}

   \subjclass[2000]{Primary: 35H10; Secondary: 35D10}

\begin{abstract}
   We consider a divergence form hypoelliptic operator
   consisting of a system of real smooth vector fields
   $X_{1},\ldots, X_{q}$
   satisfying H\"ormander condition in some
   domain $\Omega\subseteq\erren$.
   Interior $L^{p}$ estimates,
   $2\leq p<\infty$, can be obtained for weak
   solutions of the equation
   $X_j^T(a^{ij}X_iu)=X_j^T F^j,$ by assuming that
   the coefficients $a^{ij}$ belong
   locally to the space $VMO_X$ with respect to the
   Carnot--Caratheodory metric induced
   by the vector fields.
\end{abstract}

   \maketitle

   \tableofcontents

   \numberwithin{equation}{section}

\section{Introduction and main result}\label{s:introduction}
   In the present paper we obtain interior $L^{p}$ estimates for weak solutions
   of the equation $X_j^T(a^{ij}X_iu)=X_j^T F^j,$ where
   $X_{1},\ldots, X_{q}$ is a family of real vector fields
   satisfying H\"ormander's condition in some domain $\Omega\subseteq\erren,$
   and the coefficients $a^{ij}$ belong locally to the space $VMO_X,$
   with respect to the Carnot--Caratheodory metric induced
   by the vector fields. Our result generalizes, to the setting
   of hypoelliptic variational operators of H\"ormander type,
   the $L^p$ regularity results previously obtained in \cite{cfl1,difazio}.
   Indeed, in \cite{cfl1,difazio} local estimates of this kind for weak
   solutions of elliptic equation, both in divergence and non divergence form,
   are obtained by assuming that the coefficients of the operators
   belong to the space $VMO$ with respect to
   the Euclidean setting. More precisely our theorem is the following:
\begin{theorem}\label{teor:regolarita}\label{page:teo_reg}
   Let $\xuxq$ H\"ormander vector fields of step $r$ at each point of a given domain
   $\Omega\subseteq\erren,$ $q\le n$ (we can assume $n\ge 3$); moreover let
   $2\leq p<\infty.$ Let us consider the following variational equality of divergence 
   form:
\begin{equation}\label{f:lux_uguale_fcap}
   \xjt(\aij\xconi u)=\textup{div}_{X}F
\end{equation}
   where, as usual, $\textup{div}_{X}F\equiv\xjt F^j\,$ and
   moreover, $u\in W_{loc,X}^{1,p}(\Omega)$ is a weak solution of
   \eqref{f:lux_uguale_fcap} if
\begin{equation}\label{f:lux_uguale_fcap_per_esteso}
   \int_{\Omega}\aij(x)\xconi u(x)\xconj \phi(x)dx=
   \int_{\Omega}F^j(x)\xconj \phi(x)dx\quad\textup{for any
   }\phi\textup{ test in }\Omega.
\end{equation}
   \newline
   \indent
   Let us suppose that
\begin{itemize}\label{f:ipotesi_stima_finale}
   \item[$i)$] $\{\aij\}_{i,j=1,\ldots,q}$ is a symmetric measurable matrix
   defined in $\Omega$ such that $\aij \in VMO_X(B)\cap L^{\infty}(B)$
   for any open Euclidean ball $B\Subset\Omega;$
   \item[$ii)$] there exists $\nu>0$ such that
   $\frac{1}{\nu}|w|^2\leq\aij(x)w_i w_j\leq \nu
   |w|^2$ for any $w\in\erreq$ and a.e. $x \in\Omega;$
   \item[$iii)$] $F\in L_{loc}^p(\Omega,\erreq).$
\end{itemize}
   \indent
   Then, for any $\Omega'\Subset\Omega$ there exists a constant
   $c=c(\Omega',$ $\Omega,$ $\xuxq, $ $p,$ $\aij,$ $\nu)$
   and there exists an open set $\Omega'',$ $\Omega'\Subset\Omega''\Subset\Omega,$
   such that
\begin{equation}\label{f:stima_finale}
   \|u\|_{W_X^{1,p}(\Omega')}\leq c\bigg(\|u\|_{L^p(\Omega'')}+
   \|F\|_{L^p(\Omega'',\erreq)}\bigg).
\end{equation}
   \noindent
   \indent
   Note that letter $c$ denotes a generic constant that can be
   different also in the same line.
   \end{theorem}
   \noindent
   \indent
   $VMO$ functions, studied by Sarason in
   \cite{sa}, do appear first in \cite{cfl1,cfl2}, in order to obtain
   $L^p$ estimates for the solutions of uniformly elliptic equations
   in non divergence form, and later in \cite{difazio}, in non
   divergence form. In both two cases techniques rely on suitable
   representation formulas, on singular integrals depending on
   a parameter, and on their commutators with $BMO$ fuctions.
   $VMO$ condition is a type of discontinuity which implies some kind of
   average continuity: in such a sense, $VMO$ hypothesis extends classical
   theory of Agmon--Douglis--Nirenberg in \cite{agmdounir1,agmdounir2}.
   Indeed uniformly continuous bounded functions, as well the ones in $W^{1,n}$
   and $W^{1,\frac{n}{\theta}},\;\theta\in ]0,1[,$ belong to $VMO.$
   \newline
   \noindent
   \indent
   The introduction of such a family of vector fields goes back to the paper of
   H\"ormander \cite{hormander} where the author shows that hypoellipticity of
   the solution of a differential equation related to a sum of squares of vector fields
   follows from a geometric condition on the vector fields and their commutators.
   Later, Rothschild e Stein in \cite{rothstein}, deal with the problem of a natural 
   setting in which such a sum of square operators can be cast. The algebraic structures
   that do appear in this new setting are nowadays known as Carn\'ot groups;
   in particular, Euclidean spaces are a very particular cases. These are
   particular simply connected nilpotents Lie groups whose finite dimensional
   Lie algebra admits a graduated stratification in vector subspaces.
   It follows that this algebraic structure is naturally equipped with
   a family of automorphisms which generalize the standard product with
   scalars in $\erren.$ Finally, a well known theorem of Rothschild e Stein
   shows how it is possible to approximate a class of differential operators consisting 
   of a system of H\"ormander vector fields, through invariant differential
   operators defined in suitable Carn\'ot groups. From the metric viewpoint, we can 
   naturally settle these spaces in a general class of metric spaces nowadays
   known as \emph{Carn\'ot--Caratheodory Metric Spaces}, where metric is introduced
   through suitable finite families of Lipschitz vector fields.
   Such metric spaces has been intensively
   studied in the last thirty years in several setting of pure and applied mathematics
   such as degenerate elliptic differential equations, hypoelliptic differential 
   operators, sub--Riemannian manifolds, control theory, mathematical models of human 
   vision, robotics, geometric measure theory. In particular, when the vector 
   fields inducing the metric satisfy H\"ormander condition, the associate metric
   $d_{X}$ enjoys many good properties: for instance, the induced topology is actually
   the Euclidean one; all necessary properties for our purpose can be found in
   \cite{carfan}. It should be clear that, in this paper, such a metric will play a role
   becouse the coefficients $\aij$ of the operator belong to the space $VMO_{X}$ defined
   through the Carnot--Caratheodory metric induced by the vector fields associated to
   the hypoelliptic variational operator.
   \newline
   \noindent
   \indent
   Coming back to hypoelliptic operators consisting of a
   family of H\"ormander vector fields, $VMO_{X}$ functions do appear in the papers of
   Bramanti and Brandolini \cite{brmbrn1,brmbrn2}. Indeed, the coefficients
   $\aij$ are assumed to belong to the $VMO_{X}$ space with respect to the metric
   induced from the vector fields: clearly in general the space $VMO_{X}$ is different
   from the \lq\lq Euclidean\rq\rq $VMO,$ so particular metric
   proofs must be adapted in this new setting.
   Moreover, the proofs of the results in the Euclidean setting
   (see \cite{cfl1, difazio}) need several notions:
   the existence of a translation invariant fundamental solutions
   smooth away from the origin, convolution
   operators, representation formulas, parametrized singular integrals and Riesz 
   potential, commutators with $BMO$ functions, analysis on spaces of homogeneous type,
   properties of $VMO$ functions.
   In the new setting, and in particular in our case, these notions and proofs can be 
   adapted by employing the technics introduced by Rothschild e Stein, see
   \cite{brmbrn2}) so that, all in all, the properties
   of the solution of the given equation can be recovered from the
   properties of the solution of a new equation associated to a new operator defined
   locally on a suitable Carnot group, pulling back local estimates in this last setting
   to a local estimate for the solution of the given operator. In particular,
   arguing as in \cite{brmbrn2}, the use of a parametrix implies that
   coefficients $\aij$ be smooth, that is
   every function in $VMO_{X}$ should be approximated by a sequence of smooth functions:
   this is actually possible: in \cite{carfan}, in the setting of a general space of 
   homogeneous type, the space $VMO$ is defined both through balls and
   \lq\lq cubes\rq\rq, and the density property with smooth functions with respect to 
   the $BMO$ norm is proved, in the particular case of a Carn\'ot--Caratheodory spaces
   whose metric is associated to a finite family of H\"ormander vector fields. 
   \vskip15pt
   \noindent
   \indent
   The regularity result and significative properties of the space $VMO_{X}$ are 
   contained in the
   doctoral thesis discussed on December 2002. The regularity result has been announced 
   at the \emph{XVII Congresso U.M.I.} hold in
   Milan (Italy) on September 8--13, 2003, at the conference
   \emph{\lq\lq Aspetti Teorici ed Applicativi di Equazioni alle derivate parziali
   \rq\rq} hold in Maiori (Italy) on April 21--24, 2004, and at a talk given in Bologna
   (Italy) in summer 2004. The density results employed in this paper has been published 
   later in 2007, see \cite{carfan}, and so mentioned in the following.
\section{Preliminaries}\label{s:preliminaries}
   \subsection{Carn\'ot groups and Carn\'ot--Caratheodory metric spaces}
   \label{sez.:gruppi_di_Carnot}
   We refer to Section 3 of \cite{carfan} for basic definition on Carn\'ot
   groups and Carn\'ot--Caratheodory metric spaces, associated in particular to
   a family of H\"ormander vector fields; in the same section can be
   found the statement of the Ball--Box theorem, useful for our purpose. 
   \subsection{H\"ormander vector fields: theorem of Rothschild e Stein}
\label{s:campi_di_h_metr_ass}
   \indent
   Let $\xuxq$ smooth real vector fields defined on a smooth manifold.
   For $s\in\mathbb{N},$ $\,i_1,i_2,\ldots,i_{s-1},i_s$ $\in$ $\{1,2,\ldots,q\}$
   let $I=(i_1,i_2,\ldots,i_{s-1},i_s)$ and
\begin{equation}\label{f:commutat_s_esimo}
   X_I=\Big[X_{i_1},\big[X_{i_2},\ldots[X_{i_{s-1}},
   X_{i_s}]\big]\ldots\Big].
\end{equation}
   \newline
   \indent
   We say that $I$ has \emph{length} $s$ and we call
   \emph{commutator of length} $s$ any vector field such that
   $X\in \textrm{Span}\{X_I\}_{\textrm{length}\,I=s};$
   commutators of length $1$ are just the elements of the span of the vector
   fields $\xuxq.$ Suppose that for every
   $x\in M$ there exists $s(x)\in\mathbb{N}$ such that
   $\textrm{Span}\{X_I(x)\}_{\,\textrm{length}\,I\leq\,s(x)}=T_x(M);$
   then we say that the vector fields $\xuxq$
   \emph{satisfy H\"ormander condition} of step $r\in\mathbb{N}$ if
   $s(x)\le r$ for any $x\in M.$ We finally recall that the vector fields
   $\xuxq$ are \emph{free up to the step r at the point} $x\in M$ if the vectors
   ${X_I(x)}_{\,\textrm{length}\,I\leq\,r}$ are linearly independent,
   except for Jacobi's and anticommutativity relations.
   \newline
   \indent
   Let now $\gia=V_1\oplus\cdots\oplus V_r$ be the nilpotent free Lie algebra
   of step $r$ with $q=\textrm{dim}_{\erre}(V_1)$ generators,
   and let $\gig$ be the corresponding free
   Carn\'ot group (recall that,
   from a set--theoretical viewpoint, it is possible to assume
   that $\gig$ is some $\erreN$ endowed with a suitable product of Lie group
   and with a suitable family of dilations $\{\gamma_s\}_{s>0}$).
   We can denote by $\{Y_{jk^j}\}_{\substack{1\leq j\leq r\\1\leq k^j\leq n_j}}$
   a base of $V_j,$ where $j\in\{1,\ldots,r\}$ and
   $k^j\in\{1,\ldots,n_j\}$ ($n^j$ is a positive integer depending on $V_j$);
   for the sake of simplicity we shall denote by $\yuyq$ the generators
   of the first layer of $\gia;$ we can denote by $(y_{jk^j})$ the exponential 
   coordinates of the first kind of $y=\exp_{1}(Y)$ in $\gig,$ where
   $Y=\sum_{j=1}^{r}$ $\sum_{k^j=1}^{n^j}$ $y_{jk^j}$ $Y_{jk^j}$ $\in\gia$
   denotes an element of the algebra.
   Finally let $\{\delta_s\}_{s>0}$ be the family of automorphims of $\gia.$
   \newline
   \indent
   Then, the vector fields $\xuxq$ are free up to the step $r$ at the point $x\in M$
   if and only if $\textrm{dim}_{\erre}(\textrm{Span}
   \left\{X_I(x)\right\}_{\textrm{length}
   \,I\leq\,r})=$ $\textrm{dim}(\gig),$ where the last number denotes the dimension
   of $\gig$ as a smooth manifold.
   \newline
   \indent
   Let us suppose now that the vector fields do satisfy H\"ormander condition
   of step $r$ at $x_0\in M;$ let $n=\textrm{dim}(M),$
   $N=\textrm{dim}(\gig),$ $k=N-n,$
   $\widetilde{M}=M \times\mathbb{R}^k$ and let
   $\pi:\widetilde{M}\to M$ the canonical projection.
   \newline
   \indent
   Then we have the following \lq\lq lifting theorem\rq\rq of Rothschild--Stein.
\begin{theorem}\label{t:liber}
   Let $\xuxq$ smooth vector fields defined on $M$ satisfying H\"ormander condition
   of step $r$ at the point $x_0\in M.$ Then there exist
   $\{\lambda_{jl}(x,t)\}_{\substack{1\leq j\leq q\\n+1\leq l\leq N}}$
   smooth functions of the new variables
   $t_{n+1},$ $\ldots,$$t_N,$ defined
   in a neighborhood of $\xi_0=\left(x_0,0\right)$
   $\in \widetilde{U}=U\times U'\subset \widetilde{M},$ where $U$ is a
   neighborhood of $x_0$ in $M$ and $U'$ a neighborhood of $0$ in $\mathbb{R}^k,$ 
   such that the vector fields
\begin{equation*}
   \widetilde{X_j}=X_j+\sum_{l=n+1}^N\lambda_{jl}(x,t)\,\partial_{t_l},
   \quad j=1,\ldots,q,
\end{equation*}
   are free up to the step $r$ at each point of $\widetilde{U}.$$\B$
\end{theorem}
   \indent
\begin{remark}\label{o:come_riportare_la_stima_indietro}
   \textup{It is easy to verify that also the vector fields $\widetilde{X_j}$
   satisfy H\"ormander condition of step $r$ at each point
   $\xi\in\widetilde{U}$ and that it results
\begin{equation*}
   \widetilde{X_j}(f\circ \pi)=X_jf\circ \pi
\end{equation*}
   for any $f\in C^{\infty}(U)$ and for any $j=1,\ldots,q.$
   }
\end{remark}
\begin{notation}
   \emph{
   For any $f$ defined in some subset $S\subseteq\Omega,$
   we shall denote by either $f\circ\pi$ or $\widetilde{f}$
   the function defined in $S\times \errek$ that maps $\xi=(x,t)\in S\times \errek$
   to $f(x).$
   }
\end{notation}
   Let $\lambda\in\mathbb{R},$ $\lambda>0.$ A measurable function $f:\gig\to\erre$ is 
   said to be \emph{homogeneous of degree $\lambda$}
   if $f\circ\gamma_s=s^{\lambda}f$ for any $s>0;$ a differential
   operator $D$ on $\gig$ is said to be \emph{homogeneous of degree
   $\lambda$} if $D(f\circ\gamma_s)=s^{\lambda}(Df)\circ\gamma_s,$
   for any $s>0$ and for any $f\in C^{\infty}(\gig).$ Then it immediately
   follows that if $D$ and $f$ are a differential
   operator and a function, respectively homogeneous of degrees
   $\lambda_1$ e $\lambda_2,$ then $Df$ e $fD$ are a function and an operator
   homogeneous respectively of degrees $\lambda_2-\lambda_1$ and
   $\lambda_1-\lambda_2.$
   \newline
   \indent
   Let us recall now the notion of local degree at the origin (see
   \cite{rothstein} pag.272 and \cite{brmbrn2} pag.789). Let $D$ be a
   differential operator. We say that $D$ is \emph{homogeneous of local degree}
   $\leq \lambda$ if all Taylor polynomials of the coefficients of the operator
   give rise, up to a rearrangements, to a sum of differential operators of
   degree at most $\lambda$ in $\gig.$
   \newline
   \indent
   Let us suppose now $\xuxqtilde$ be free vector fields up to 
   the step $r$ at a point $\xi_0$ of a smooth manifold $\widetilde{M}.$
   Then $\textrm{dim}_{\erre}(\textrm{Span}
   \{\widetilde{X_I}(\xi_0)\}_{\textrm{length}\,I\leq\,r})=
   \textrm{dim}(\gig).$ Now, if all $\xjtilde$ were invariant on $\gig,$
   it would be possible to identify them with the elements of (the first layer of)
   $\gia$ and, consequently, recover the elements of $\gig$ through the before mentioned
   exponential coordinates; this is not possible in general because the vector fields 
   are not invariant in general; nevertheless, if we choose
   $\{\widetilde{X_{jk^j}}\}$ such that $\textrm{Span}
   \{\widetilde{X_I}(\xi_0)\}_{\textrm{length}\,I\leq\,r}=
   T_{\xi_0}(\widetilde{M}),$
   we can consider the mapping that, for any
   $N$--tuple of real numbers $y=(y_{jk^j})$ in a fixed closed ball $B$
   around $0$ sufficiently small, maps any $\xi$ in a compact neighborhood
   $\widetilde{U'}$ of $\xi_0$ to the element
   $\eta=\exp_{\widetilde{X}}(1),\;
   \widetilde{X}=\sum_{j=1}^{r}\sum_{k^j=1}^{n^j}y_{jk^j}\widetilde{X_{jk^j}}:$
   indeed, under these conditions, we have the existence, unicity and
   $C^{\infty}$ dependence from $\xi$ ed $y$ of the Cauchy problem
   $\phi'(t)=\widetilde{X}(\phi(t))$
   and $\phi(0)=\xi,$ for $t$ in a sufficiently small neighborhood of $0,$ and
   for any fixed $\xi\in \widetilde{U'},\;y\in B.$
   Then we have
   $\phi(t)=\exp
   \left(t\big(\sum_{j=1}^{r}\sum_{k^j=1}^{n^j}y_{jk^j}
   \widetilde{X_{jk^j}}\big)\right)\xi;$
   in particular $\eta=\phi(1)=\exp
   \left(\sum_{j=1}^{r}\sum_{k^j=1}^{n^j}y_{jk^j}
   \widetilde{X_{jk^j}}\right)\xi$
   and
\begin{equation*}
   y=\Theta_\xi(\eta)\quad\textrm{e}\quad
   \eta=\Theta_\xi^{-1}(y),\quad\textrm{for any }\xi,\eta\in
   \widetilde{U'},y\in B.
\end{equation*}
   \indent
   The mapping $\Theta_\xi$ then behaves like a
   coordinate map; indeed through
   $\Theta_\xi$ we can think of $\xjtilde$ as defined on $\gig$ and
   consequently to approximate them with the left invariant vector fields $\yj$
   on $\gig$ (which, in their turn, can be chosen so that they agree with
   the $j$--th partial derivatives at the origin (see \cite{rothstein} pag.
   272; see also \cite{sc} for the following formulation); we can also think of the
   $y=(y_{jk^j})$ as a system of canonical coordinates, depending only on the vector 
   fields $\{\widetilde{X_{jk^j}}\}_{\substack{1\leq j\leq r\\1\leq k^j\leq n_j}}.$
   \newline
   \indent
   Then we have the following approximation theorem of Rothschild and Stein.
\begin{theorem}\label{t:appr}
   Let $\xuxqtilde$ real smooth vector fields defined on a smooth manifold
   $\widetilde{M}$ and let $\xi_0\in\widetilde{M}.$ Let us suppose that the vector 
   fields $\xuxqtilde$ do satisfy the H\"ormander hypothesis of step $r;$
   moreover we assume that the vector fields are free of step $r$ at the same point
   $\xi_0.$ Let us choose the vector fields
   $\{\widetilde{X_{jk^j}}\}_{\substack{1\leq j\leq r\\1\leq k^j\leq n_j}}$
   as before and let us denote by $(y_{jk^j})$ the associated system of canonical 
   coordinates. Let $\gia$ and $\gig$ respectively the free Lie algebra with $q$
   generators of step $r$ and the associated free Carn\'ot group.
   \newline
   \indent
   Then it is possible to choose a base
   $\{Y_{jk^j}\}_{\substack{1\leq j\leq r\\1\leq k^j\leq n_j}}$ of
   $\gia$ such that $Y_j(0)=\frac{\partial}{\partial y_j},
   \;j=1,\ldots,q,$ a neighborhood $U$ of $0$ in $\gig,$
   two open neighborhood $W,V$ of $\xi_0$ in $\widetilde{M},$ $W\Subset V,$
   such that the following facts hold:
\begin{itemize}
   \item[i)] $\Theta_\xi\mid V$ is a diffeomorphism between $V$ and
   $\Theta_\xi(V)$ for any $\xi\in V;$
   \item[ii)] $\Theta_\xi(V)\supseteq U$ for any $\xi\in W;$
   \item[iii)] the mapping $\Theta: V\times V\rightarrow\gig$
   defined by the position $\Theta(\xi,\eta)=\Theta_\xi(\eta)$
   belongs to $C^\infty\left(V\times V\right);$
   \item[iv)] for any fixed $\xi\in W$ the mapping
   $\eta\to\Theta_\xi(\eta)=\Theta(\xi,\eta)=(y_{jk^j}),\;\eta\in W,$
   is a coordinate map for $W$ and
   $(\Theta_\xi)_\ast\widetilde{X_i}=Y_i+R^\xi_i$ in $U,$ dove $R^\xi_i$
   is a real smooth vector fields of local degree $\leq 0,$
   with $C^{\infty}$ dependence on $\xi\in W;$
   more precisely it means that for any
   $\xi\in W$ and for any $f\in C^\infty(\gig)$ it results
\begin{equation}
   \widetilde{X_i}\circ(f\circ\Theta_\xi)=
   (Y_i\circ f)\circ \Theta_\xi+(R^\xi_i\circ f)\circ\Theta_\xi.
\end{equation}
\end{itemize}
   In general, for any couple of indexes $j$ e $k^j$ and,
\begin{equation*}
   (\Theta_\xi)_\ast\widetilde{X_{jk^j}}=Y_{jk^j}+R^\xi_{jk^j}
\end{equation*}
   where $R^\xi_{jk^j}$ is a real smooth vector field of local degree
   $\leq j-1$ with $C^{\infty}$ dependence on $\xi\in W.$$\B$
\end{theorem}
   Recall that $(\Theta_\xi)_\ast$ is the mapping induced by $\Theta_\xi$
   on the fiber bundle and defined by the position
   $((\Theta_\xi)_\ast\widetilde{X})f=(\widetilde{X}\circ(f\circ\Theta_\xi))
   \circ\Theta_\xi^{-1},$ for any vector field and for any $f\in C^{\infty}(V).$
\subsection{Introduction of a quasi-metric equivalent to C--C metric}\label{ss:intr_metr_ro}
   In what follows we set $M=\Omega$ and
   $\widetilde{M}=\Omega\times\errek\subseteq\erreN,$ then our neighborhoods are
   $C\textrm{--}\,C$ balls which are open sets
   in any one of the topologies $\tau_{Euclidean}$
   and $\tau_{C\textrm{--}\,C}.$ We shall denote by $Q$ the homogeneous dimension
   as a doubling spaces in $\gig.$ It will be more useful to introduce a
   quasi-metric $\ro$ in $V,$ equivalent to $d_{\widetilde{X}}.$
   Indeed, referring the reader to \cite{sc} and \cite{brmbrn2}, let us recall
   main definitions.
\begin{theorem}\label{t:teta}
   Let $V$ and $W$ neighborhood of $\xi_0$ as in Theorem \eqref{t:appr}.
   If we set
\begin{equation*}
   \rho\left(\xi,\eta\right)=\left\|\Theta(\xi,\eta)\right\|_{\gig}
   \quad
   for\;any\;\xi,\eta\in V
\end{equation*}
   then the following properties hold
\begin{itemize}
   \item [i)] $\Theta(\xi,\eta)=\Theta(\eta,\xi)^{-1}=-\Theta(\eta,\xi);$
   \item [ii)] $\rho\left(\zeta,\eta\right)\leq
   c\,\big(\rho\left(\xi,\zeta\right)+\rho\left(\eta,\xi\right)\big)\quad$
   for any $\xi,\eta\in V$ such that $\rho(\xi,\eta)\leq 1$ and
   $\rho(\xi,\zeta)\leq 1;$
   \item [iii)] there exist four positive smooth functions
   $V\ni\xi\to \zeta(\xi),\omega(\xi),$
   $V\ni\eta\to h(\eta),$ $U\ni y\to j(y)$
   and a positive constant $c$ such that
   $1/c\leq \zeta,\omega,h,j\leq c$ respectively on $V$ the first three ones, on
   $U$ the last one, and
   $J_{\xi}(\eta)=\zeta(\xi)h(\eta)$ e
   $J^{-1}_{\xi}(y)=\omega(\xi)j(y),$
   where $J_{\xi}(\eta)$ and $J^{-1}_{\xi}(y)$ are respectively
   the jacobians of the mappings $y=\Theta_\xi(\eta)$
   and $\eta=\Theta_\xi^{-1}(y).$$\B$
\end{itemize}
\end{theorem}
   Let now $x_0\in \Omega,$
   and $\xi_0=(x_0,0)\in \widetilde{\Omega};$
   we can assume that $V=B\times R$ where $B$ is an open Euclidean ball around
   $x_0$ and $I$ is an open rectangle in $\errek;$ 
   consequently we can consider the following three quasi--metric spaces
   with respective Lebesgue measures
\begin{equation}\label{f:tre_sp_dou}
   (B,d_{X},dx),\quad (V,d_{\widetilde{X}},d\xi),\quad (V,\ro,d\xi).
\end{equation}
   Then from Theorem \eqref{t:teta} we have the following proposition.
\begin{proposition}
   According to above notation, for any
   $\xi_0=(x_0,0)\in\widetilde{\Omega}$,
   there exists an Euclidean neighborhood $V$ of the type
   $B\times I,$ where $B$ is an Euclidean open ball around $x_0$
   and $I$ is an open rectangle in $\errek$ around $0,$
   such that
   $(V,\ro,d\xi)$ is a bounded doubling space. In particular, Lebesgue measure
   $d\xi\equiv m_N$ is a doubling measure and, for any open ball $B_{\ro}\subset V$
   related to the metric $\ro$ we have $m_N(B_{\ro})\approx r^Q.$$\B$
\end{proposition}
   Finally, as in Lemma 7 of \cite{sc} we have the following proposition
   (see \cite{nastwa2} for the proof).
\begin{proposition}
   According to above notation, there exists $c>0$ such that
\begin{equation*}
   \frac{1}{c}\ro(\eta,\xi)\leq d_{\widetilde{X}} (\eta,\xi)\leq
   c\ro(\eta,\xi)\quad per\,ogni\;\eta,\xi\in V.\B
\end{equation*}
\end{proposition}
   So from the metric viewpoint, the metric spaces
   $(V,d_{\widetilde{X}},d\xi)$ and $(V,\ro,d\xi)$ are equivalent.
   Let us now compare the spaces $BMO_X(B)$ (resp. $VMO_X(B)$) defined in $B$
   with respect to the metric induced from the vector fields $\xuxq,$
   endowed with the corresponding Lebesgue measure, with the spaces
   $BMO_{\xtilde}(V)$ (resp. $VMO_{\xtilde}(V)$)
   of functions defined on $B$ with respect to the metric induced
   by the vector fields $\xuxqtilde,$ and endowed with the correpsonding Lebesgue 
   measure. The arguments are taken from pagg. 793--794 of \cite{brmbrn2},
   then we just recall the following propositions.
\begin{proposition}\label{p:confr_mis_palle_sancal}
   According to above notation, for any $x\in\Omega$ and
   $r>0$ such that $B_{d_{X}}(x,r)\subset B$ and
   $B_{d_{\widetilde{X}}}(\xi,r)\subset V,$
   if $\pi:\widetilde{\Omega}\to\Omega$ and
   $\pi_k:\widetilde{\Omega}\to\errek$ are the canonical projections,
   denoting by $m_n,\,m_k,\,m_N$ respective Lebesgue measures,
   and setting $C_k=m_k(\pi_k(B_{d_{\widetilde{X}}}(\xi,r)))$,
   it results
\begin{itemize}
   \item $\pi(B_{d_{\widetilde{X}}}(\xi,r))=B_{d_{X}}(x,r);$
   \item $d_{\widetilde{X}}(\xi,\xi ')\geq d_{X}(x,x')\quad
   \forall x,x'\in \Omega;$
   \item there exists $C>0$ such that, for $r$ small enough,
   $\frac{1}{CC_k}m_N(B_{d_{\widetilde{X}}}(\xi,r))\leq
   m_n(B_{d_{X}}(x,r))\leq\frac{C}{C_k}m_N(B_{d_{\widetilde{X}}}(\xi,r)).$$\B$
\end{itemize}
\end{proposition}
\begin{proposition}\label{p:equiv_vmo}
   According to above notation, if $f:\Omega\to\erre$
   is a measurable function, then $f\in BMO_X(B)$ (resp. $f\in VMO_X(B)$)
   with respect to metric $d_{X}$ and Lebesgue measure $m_n$ if and only if
   $f\circ\pi\in BMO_{\xtilde}(V)$ (resp. $f\circ\pi\in VMO_{\xtilde}(V)$)
   with respect to metric $\ro$ Lebesgue measure $m_N.$$\B$
\end{proposition}
\subsection{Differential operators, fundamental solutions and
   parametrices}\label{ss:op_diff_sol_fond_par}
   In this section our space of homogeneous type will be
   $(V,d_{\widetilde{X}},d\xi)\equiv(V,\ro,d\xi),$ according to notations of Section
   \eqref{ss:intr_metr_ro}. We recall results of Sections 2.1. e 3.2. in
   \cite{brmbrn2}. Let $L$ a given differential operator consisting of a family
   $X_0,\xuxq$ of H\"ormander vector fields defined on a given open set
   $\Omega\subset\erren.$ For instance let either $L=\sum_{i=1}^n X_i^2+X_0,$
   or $L=\aij X_i X_j,$ where the coefficients belong to $C^\infty(\Omega).$
   Arguing as in \cite{rothstein} we will recover properties of the operator
   $L$ from the properties of a new operator $\widetilde{L},$ consisting
   of the $\widetilde{X_0},\xuxqtilde;$ this last operator, in its turn,
   has much more properties useful for our purpose becouse it can be written
   as a sum of two more operators: the first one consists of left invariant
   vector fields defined on a suitable nilpotent group $\gigno$,
   the second one, defined in $\gigno,$  is a smooth operator
   of local degree equal to zero so that results of Folland
   (see Teorema 2.1, Section 2 in \cite{fo1}) can be applied, for the
   existence of a fundamental solution. Finally, quoting Christ
   (see \cite{christ2}, Example 8, Pag.96) we need some more analysis
   to obtain our estimates; more precisely we need to construct,
   through a suitable coordinate map (see Theorem \eqref{t:appr}), two
   parametrices (see \S15. in \cite{rothstein}) that behave, in our
   case, as left and right partial inverse of the operator $\widetilde{L}$,
   up to a finite number of operators (depending only on $\widetilde{L}$) which,
   in their turn, are the analogous of classical integral with either singular
   kernel, or fractional or Riesz potential. So, given the estimates with new vector 
   fields $\widetilde{X_0},\xuxqtilde,$ we can recover the original estimate
   (see Remark \eqref{o:come_riportare_la_stima_indietro}).
   In our case, the operator is $L=\lux$ where, as before, $\xuxq$ is a family
   of H\"ormander vector fields defined on a given open set $\Omega\subset\erren,$
   whose coefficients belong locally to the class $VMO_X$ with respect to
   $C\textrm{--}\,C$ metric induced from the vector fields $\xuxq;$
   we consider weak solutions for the equation. Then,
   through vector fields $\xuxqtilde,$ we pass from the equation
   associated to the operator $\lux$ to the one associated to the
   divergence form operator $\luxtilde$ consisting of the vector fields
   $\xuxqtilde$ (note that, in this case, thanks to Proposition \eqref{p:equiv_vmo},
   coefficients $\aijtilde$ belong locally to $VMO_{\xtilde}$
   with respect to $C\textrm{--}\,C$ metric induced by vector fields $\xuxqtilde;$
   if we are able to obtain local estimate for the solution of the new equation,
   we can obtain the requested estimates. So we need estimates for the
   solutions of the equation associated to the operator
   $\widetilde{L}=\luxtilde$ with coefficients locally in $VMO_{\xtilde}$
   with respect to $C\textrm{--}\,C$ metric associated to vector fields $\xuxqtilde.$
   Now, thanks to results of Section \eqref{sez.:sp_om_e_dens}, it is possible
   to approximate locallly coefficients $\aijtilde$ with smooth functions; then
   the divergence form operator agrees, up to some low order terms which belong to
   the span of the vector fields $\xuxqtilde,$ with the non divergence form 
   hypoelliptic one (see Teorema 1.11 in \cite{brmbrn1}). For the main part of
   this operator we can consider the parametrix adapted by the authors 
   in \cite{brmbrn2} from the original one of Rothschild--Stein in \cite{rothstein}, 
   pag. 296. Then, by assuming that coefficients $\aijtilde$ are smooth, we can
   obtain a first estimate for a test solution $u:$ this is the argument
   of \emph{Step 2} of pag. \pageref{page:passo 2}.
   \newline
   \indent
   Let us now recall the results in the form useful for this purpose. 
\begin{theorem}\label{t:teor_folland_sol_fond}
   Let $\xuxq$ be H\"ormander vector fields, let $\xuxqtilde$
   be the free vector fields associated and let $Y_1,\ldots,Y_q$ $\in$ $\gia$
   the approximating left invariant vector fields,
   according to Theorem \eqref{t:appr}.
   Then, for any fixed $\xi_0\in\gig,$ the operator
   $\ajitilde(\xi_0)Y_jY_i$ is hypoelliptic jointly with its transposed.
   Under this conditions there exists $\Gamma_0\equiv\Gamma_{\xi_0}\in
   C^\infty(\gig\setminus\{0\}),$ homogeneous in $\gig$ of degree $2-Q$
   and such that for any test funciotn $\phi$ in $\gig$ and any $\xi\in\gig$
   it results
\begin{equation*}
   \phi(\xi)=\int_{\gig}
   \Gamma_0(\eta^{-1}\xi)(\ajitilde(\xi_0)Y_jY_i\phi)(\eta)d\eta.
\end{equation*}
   Moreover, for any $i,j=1,\ldots,q$ there exists constants
   $\alpha_{ij}(\xi_0)$ such that for any $\xi\in\gig$ we have
\begin{equation*}
   Y_jY_i\phi(\xi)=\textup{P.V.}
   \int_{\gig}Y_jY_i\Gamma_0(\eta^{-1}\xi)(\ajitilde(\xi_0)Y_jY_i\phi)(\eta)d\eta
   +\alpha_{ji}(\xi_0)(\ajitilde(\xi_0)Y_jY_i\phi)(\xi),
\end{equation*}
   $Y_jY_i\Gamma_0\in C^\infty(\gig\setminus\{0\}),$ where
   $Y_jY_i\Gamma_0$ is homogeneous in $\gig$ of degree $-Q.$
   Moreover $\sup_{\xi\in\erreN}|\alpha_{ji}|<\infty.$$\B$
\end{theorem}
   \indent
   According to notations in previous sections let us recall definitions
   and main properties of \lq\lq operators of type $0,1,2$\rq\rq.
\begin{definition}[Kernel and operators of type $0,1,2$]\label{d:oper_tipo_elle}
   Let $\xi_0\in V$ be fixed.
   We say that $K_{0,0}(\xi,\eta)\equiv K_{\xi_0,0}(\xi,\eta)$
   (resp. $K_{0,1}(\xi,\eta)\equiv K_{\xi_0,1}(\xi,\eta),
   \,K_{0,2}(\xi,\eta)$ $\equiv$ $K_{\xi_0,2}(\xi,\eta)$)
   is a \emph{frozen kernel at}
   $\xi_0$ \emph{of Type} $0$ (resp. $1,2$) if, according to notations
   of Theorem \eqref{t:teor_folland_sol_fond}, for any $m\in\mathbb{N},$
   it can be written as a finite sum of the kind
\begin{equation*}
   \Big[a_0(\xi)(D_0\Gamma_0)(\tetaxe)b_0(\eta)\Big]+
\end{equation*}
\begin{equation*}
   +\Big[a_1(\xi)(D_1\Gamma_0)(\tetaex)b_1(\eta)\Big]+
   \cdots+\Big[a_s(\xi)(D_s\Gamma_0)(\tetaex)b_s(\eta)\Big]
\end{equation*}
   with $s=s(m)\in\mathbb{N},$ $a_i,b_i$ test functions in $V$
   for any $i=0,1,\ldots,s,$
   $D_1,\cdots,D_s$ differential operators homogeneous of degree less or equal than
   2 (resp. 1,0), and $D_0$ is differential operator such that $D_0\Gamma_0\in C^m(V).$
   \newline
   \indent
   Let now $\phi\in C_c^\infty(V).$
   \newline
   \indent
   We say that $T_{0,0}\equiv T_{\xi_0,0}$
   is a \emph{frozen operator at}
   $\xi_0$ \emph{of Type} $0$
   if there exists a bounded measurable function
   $\alpha_0\equiv\alpha_{\xi_0}$ such that, for any
   $\xi\in V,$
\begin{equation*}
   (T_{0,0}\phi)(\xi)=\textup{P.V.}\int_V
   K_{0,0}(\xi,\eta)\phi(\eta)d\eta+\alpha_0(\xi)\phi(\xi);
\end{equation*}
   We say then that $T_{0,1}\equiv T_{\xi_0,1}$
   (resp. $T_{0,2}\equiv T_{\xi_0,2}$)
   is a \emph{frozen operator at}
   $\xi_0$ \emph{of Type} $1$ (resp. $2$) if
\begin{equation*}
   (T_{0,1}\phi)(\xi)=\int_V K_{0,1}(\xi,\eta)\phi(\eta)d\eta
   \quad
   \left(
   resp.\;(T_{0,2}\phi)(\xi)=\int_V K_{0,2}(\xi,\eta)\phi(\eta)d\eta
   \right).
\end{equation*}
   If, for any $k=0,1,2,$ $K_{0,k}(\xi,\eta)
   \equiv K_{\xi_0,k}(\xi,\eta)$ is a frozen kernel at
   $\xi_0$ of type $k,$ then we say that
   $K_k(\xi,\eta)\equiv K_{\xi,k}(\xi,\eta)$ is \emph{a Kernel of Type}
   $k.$
   \newline
   \indent
   Finally we say that $T_0\equiv T_{\xi,0}$
   is an \emph{Operator of Type} $0$ if there exists a bounded measurable function
   $\alpha_0(\xi)\equiv\alpha_{\xi_0}(\xi)$ such that, setting
   $\alpha(\xi)\equiv\alpha_{\xi}(\xi),$ it results
\begin{equation*}
   (T_0\phi)(\xi)=\textup{P.V.}\int_V
   K_0(\xi,\eta)\phi(\eta)d\eta+\alpha(\xi)\phi(\xi);
\end{equation*}
   analogously we say that $T_1\equiv T_{\xi,1}$
   (resp. $T_2\equiv T_{\xi,2}$)
   is an \emph{Operator of Type} $1$ (risp. $2$) if
\begin{equation*}
   (T_1\phi)(\xi)=\int_V K_1(\xi,\eta)\phi(\eta)d\eta
   \quad
   \left(
   resp.\;(T_2\phi)(\xi)=\int_V K_2(\xi,\eta)\phi(\eta)d\eta
   \right).
\end{equation*}
\end{definition}
   According with above notations, the following facts hold.
\begin{lemma}\label{l:tipo_i_tipo_i-1}
   If, for any $k=1,2,$ $K_{0,k}(\xi,\eta)$ is a frozen kernel at
   $\xi_0$ of type $k,$ then
   $(\xitilde K_{0,k}(\cdot,\eta))(\xi)$ ia a frozen kernel at $\xi_0$
   of type $k-1.$
   \newline
   \indent
   If, for any $k=1,2,$ $T_{0,k}$ is a frozen kernel at $\xi_0$ of type $k,$
   then $\xitilde T_{0,k}$ is a frozen operator at $\xi_0$ of type $k-1.$$\B$
\end{lemma}
\begin{example}\label{e:vari_tipi_di_nuclei}
   \textup{
   We recall that, for instance, fixed
   $\xi_0\in V,$ if $i,j=1,\ldots,q,$ then
\begin{itemize}
   \item[--] $a(\xi)\Gamma_0(\tetaex)b(\eta),\;$
   \item[--] $a(\xi)(R_i^\eta\Gamma_0)(\tetaex)b(\eta),$
\end{itemize}
   are frozen kernel at $\xi_0$ of type $2;$ while,
\begin{itemize}
   \item[--] $a(\xi)(Y_i\Gamma_0)(\tetaex)b(\eta),$
   \item[--] $a(\xi)(Y_iR_j^\eta\Gamma_0)(\tetaex)b(\eta),$
   \item[--] $a(\xi)(R_i^\eta R_j^\eta\Gamma_0)(\tetaex)b(\eta),$
\end{itemize}
   are frozen kernel at $\xi_0$ of type $1.$
   }
\end{example}
   Let us conclude the present section with the following theorems
   whose proofs is either taken or adapted from the ones in \cite{brmbrn2}).
\begin{theorem}\label{t:xtilde_con_operatore_scambiano}
   Let $T_{0,k}$ be a frozen operator at $\xi_0\in V$ of type $k=0,1,2.$
   Then, for any vector field $\xitilde$ there exist $q+1$ operators
   $T^0_{0,k},T^1_{0,k},\ldots,T^q_{0,k}$
   frozen at $\xi_0,$ of the same type $k$ of $T_{0,k}$ such that
\begin{equation*}
   \xitilde T_{0,k}=\sum_{h=1}^q T^h_{0,k}\widetilde{X_h}+T^0_{0,k}
\end{equation*}
\end{theorem}
\begin{theorem}\label{t:limitatez_int_sing}
   Let $T_0$ an operator of type $0$ and $1<p<\infty.$ Then there exists a constant
   $c\equiv c(T_0,p)$ such that for any $u\in L^p(V)$ and for any
   $a\in BMO_{\xtilde}(V)$ it results
\begin{enumerate}
   \item $\|T_0u\|_{L^p(V)}\leq c\|u\|_{L^p(V)};$
   \item $\|[T_0,a](u)\|_{L^p(V)}
   \leq c\|a\|_{BMO_{\xtilde}(V)}\|u\|_{L^p(V)};$
   \item if moreover $a\in VMO_{\xtilde}(V)$ and $\ep>0,$ then there
   exists $r>0$ depending on $p,T_0,\ep$ and $a$ such that,
   for any ball $\widetilde{B}$ associated to the metric induced from the
   vector fields $\xuxqtilde,$ if $\textup{supp }u\subset\widetilde{B}$ then
\begin{equation*}
   \|[T_0,a](u)\|_{L^p(V)}\leq \ep\|u\|_{L^p(V)},
\end{equation*}
\end{enumerate}
   where $BMO_{\xtilde}(V)$ and $VMO_{\xtilde}(V)$
   denote function spaces as in Section \eqref{sez.:sp_om_e_dens},
   with respect to $C\textrm{--}\,C$ metric induced from vector fields
   $\xuxqtilde,$ and $[T_0,a]$ denotes the commutator which maps $u\in L^p(V)$
   into $T(au)-aT(u).$$\B$
\end{theorem}
\begin{theorem}
   Let $T_k$ be an operator of type $k=1,2$ and $1<p<\frac{Q}{k}.$ Then there exixts
   a constant $c\equiv c(T_k,p)$ such that, if
   $\frac{1}{q}=\frac{1}{p}+\frac{k}{Q},$
   then, for any $u\in L^p(V)$ it results
\begin{equation*}
   \|T_ku\|_{L^p(V)}\leq c\|u\|_{L^q(V)}.\B
\end{equation*}
\end{theorem}
\begin{notation}\label{notaz:p_stars}
   \textup{
   Set $\frac{1}{p_\ast}=\frac{1}{p}+\frac{1}{Q}$ and, for
   $p<Q,$ $\frac{1}{p^{\ast}}=\frac{1}{p}-\frac{1}{Q}$
   (see \emph{Step 4} pag. \eqref{s:regularity}). It is easy verified that
   $p_\ast<p<p^\ast,$ $(p_\ast)^{\ast}=(p^\ast)_\ast=p$ and,
   moreover, the mapping that associate to $p$ any one of the corresponding
   \lq\lq star $p$\rq\rq is order preserving in the real numbers;
   for instance, from $p_\ast<p<p^\ast$ it follows that $p<p^\ast<p^{\ast\ast}$
   then $p^\ast<p^{\ast\ast}<p^{\ast\ast\ast},$ and so on
   }
\end{notation}
\begin{corollary}\label{cor:i_due_tipi_vanno_a_star_basso}
   Let $T_k$ un operator of type $k=1,2$ and $1<p<\frac{Q}{2}.$
   Then there exists a constant $c\equiv c(T_1,T_2,p)$ and $V$
   such that for any $u\in L^p(V)$ it results
\begin{align*}
   &\|T_1u\|_{L^p(V)}\leq c\|u\|_{L^{p_{_{\ast}}}(V)}\\
   &\|T_2u\|_{L^p(V)}\leq c\|u\|_{L^{p_{_{\ast}}}(V)}.\B
\end{align*}
\end{corollary}
\subsection{Sobolev spaces associated to a family of Lipschitz vector fields}\label{ss:sp_sob_cc}
   Let $\Omega\subset\erren$ an open set and
   $Y:\Omega\to\erren$ a Lipschitz vector field
\begin{equation*}
   Y(x)=\sum_{i=1}^{n}b_{i}(x)\partial_{i}\equiv(b_{1}(x),\ldots,b_{n}(x))
   \quad \forall x\in \Omega.
\end{equation*}
   Assuming, for instance, Lipschitz regularity for $\partial\Omega,$ if
   $f,b_i\in C^1(\overline{\Omega})\;,i=1,\ldots,n,$ and
   $\phi \in C_c^\infty(\Omega)$ then
\begin{equation*}
   \int_\Omega Yf\phi\,dx
   =\sum_{i=1}^n\left[\int_\Omega \partial_i(b_if\phi)dx
   - \int_\Omega f\partial_i(b_i\phi)dx\right]
   =\int_\Omega f\left(-\sum_{i=1}^n \partial_i(b_i\phi)\right)dx.
\end{equation*}
   If we set $Y^T=-\sum_{i=1}^n \partial_i(b_i\,\cdot)$
   then we can write
\begin{equation*}
   \int_\Omega Yf\phi\,dx = \int_\Omega fY^T\phi\,dx
\end{equation*}
   so it suffices to request for any $b_i$ to be locally Lipschitz.
   If $Y=\sum_{i=1}^nb_{i}\partial_{i}$ is a locally Lipschits vector fields
   on $\Omega$ and $f,\,g\in L_{loc}^1(\Omega),$ we say that $g$ is the
   partial derivative along $Y,$ and we write $Yf=g,$ if
\begin{equation*}
   \int_\Omega  g\phi\,dx=\int_\Omega f Y^T\phi\,dx
\end{equation*}
   for any $\phi \in C_c^\infty(\Omega).$
   Let $X=(X_1,\ldots,X_q)$ a family of locally Lipschitz vector fields
   on $\Omega.$ if $f\in L^1_{loc}(\Omega)$ has partial derivatives along
   $X_j\;\forall j=1,\ldots,q,$ let us denote by
   $Xf=(X_1f,\ldots,X_qf)$ the \emph{weak gradient of $f$}. Moreover we set
   $|Xf|=\left(|X_1f|^2+\cdots+|X_qf|^2\right)^{\frac{1}{2}}.$
   \newline
   \noindent
   \indent
   Let $1\leq p<\infty,$ and let $\{$$X_1,\ldots,X_q$$\}$ be a family of locally Lipschitz
   vector fields on $\Omega.$
\begin{definition}\label{d:spaz_sob_con_campi}
   The Sobolev space $W_{X}^{1,p}(\Omega)$
   is the space of all function $f:\Omega\to\mathbb{R}$
   such that $f\in L^p(\Omega)$ and, for any $j=1,\ldots,q,$
   $X_jf$ do exists in the weak sense and belong to $L^p(\Omega).$
\end{definition}
   $W_{X}^{1,p}(\Omega)$ is a Banach space with the norm
   $$\|f\|_{W_X^{1,p}(\Omega)}=
   \Big(\|f\|_{L^p(\Omega)}^p+
   \sum_{j=1}^{q}\|X_jf\|_{L^p(\Omega)}^p\Big)^{\frac{1}{p}}\,,$$
   for any $1\leq p<\infty,$ or equivalently with the norm
   $\|f\|_{L^p(\Omega)}+$ $\sum_{j=1}^{q}$ $\|X_jf\|_{L^p(\Omega)}.$
   Analogous definitions hold for the local Sobolev spaces $W_{loc,X}^{1,p}(\Omega),$
   and for the subspace of all functions zero on $\partial\Omega.$ Let us denote
   by $\overset{\circ}{W}_{X}^{1,p}(\Omega)$ the closure of
   $C_c^\infty(\Omega)$ in $W_{X}^{1,p}(\Omega).$
   The following proposition recall basic properties of Sobolev spaces.
\begin{proposition}\label{p:prop der deb con campi}
   Let $f,\,g \in W_{X}^{1,p}(\Omega).$ Then
\begin{enumerate}
   \item[(i)] for any $\lambda$ e $\mu\in\mathbb{R},\;\lambda f +\mu g
   \in W_{X}^{1,p}(\Omega)$ e $X_j(\lambda f +\mu g)=\lambda
   X_jf+\mu X_jg\quad\forall j=1,\ldots,q\,;$
   \item[(ii)] If $U$ is an open subset of $\Omega$ then
   $f\in W_{X}^{1,p}(U)\,;$
   \item[(iii)] If $\zeta\in C_c^\infty(\Omega)$ then $\zeta f\in
   \overset{\circ}{W}_{X}^{1,p}(\Omega)$ e $X_j(\zeta f)
   =\zeta X_jf+f X_j\zeta\quad\forall j=1,\ldots,q\,.$
\end{enumerate}
\end{proposition}
   \noindent
   \indent
   Let us conclude this section recalling that for $p=2,$
   $W_{X}^{1,2}(\Omega)\equiv H_X^1(\Omega)$ and
   $\overset{\circ}{W}_{X}^{1,2}(\Omega)\equiv\overset{\circ}{H}_X^1(\Omega)$
   have a natural structure of Hilbert space.
   In this setting, the weak formulation of Dirichlet problem
\begin{equation*}
   {
\begin{cases}
   \lux=f \\
   u \in\overset{\circ}{H}_X^1(\Omega)
\end{cases}}
\end{equation*}
   with $f\in H_X^1(\Omega),$ $\aij$ are bounded and measurable functions,
   and by assuming uniform ellipticity,
   has an unique solution: the proof follows standard
   Lax--Milgram lemma; analogously, classic $L^2$ regularity theory holds.
\subsection{Space of homogeneous type: the space $BMO_{X}$ and $VMO_{X}$}\label{sez.:sp_om_e_dens}
   For these definitions and results we refer the
   reader to the whole paper \cite{carfan}.
\section{Regularity Result}\label{s:regularity}
   In this section we give a detailed sketch of the proof of Theorem
   \eqref{teor:regolarita} of pag. \pageref{page:teo_reg},
   postponing at the end main calculations.
   \proof[Sketch of the proof]
   According to notation of Theorem \eqref{t:liber}, let $W=B\times I$
   where $B=B(x,r)$ is an Euclidean open ball, $x\in\Omega',$
   $I$ is an open rectangle in $\mathbb{R}^k$ centered in
   zero, and $r>0$ is sufficiently small. It suffices to verify that for any
   $\xi\in W$ there exist two open $C\textrm{--}\,C$ balls
   $B'_{\xtilde}\Subset B''_{\xtilde}\Subset W$ around $\xi$ such that
\begin{equation}\label{f:stima_finale_palla_palla_doppia}
   \|\utilde\|_{W_{\xtilde}^{1,p}(B'_{\xtilde})}
   \leq c\Big(\|\utilde\|_{L^p(B''_{\xtilde})}+
   \|\fcaptilde\|_{L^p(B''_{\xtilde},\erreq)}\Big).
\end{equation}
\begin{itemize}
   \item[\small{(\emph{Step 1})}] Let $p\geq 2$ fixed.
   The function $u\in W_{X}^{1,p}(\Omega)$ is a weak solution
   in $\Omega$ of the equation \eqref{f:lux_uguale_fcap}; then
   $\utilde\in W_{\xtilde}^{1,p}(W)$ is a weak solution in $W$ of the equation
\begin{equation}\label{f:luxtilde_uguale_fcaptilde_piu_f}
   \xjttilde(\aijtilde\xitilde\utilde)=\xjttilde\fcapjtilde+f,
\end{equation}
   \indent
   where $\xuxqtilde$ are the free vector fields as in Theorem \eqref{t:liber},
\begin{equation}\label{f:f}
   f=(\aij\xitilde\utilde-\fcaptilde^j)g_j,
\end{equation}
   and $g_j$ are smooth functions in $W,$ defined through the positions
\begin{equation}\label{f:g_con_j}
   g_j=-\Big(\paj\lambda_{j n+1}+\cdots+\paj\lambda_{j N}\Big)=
   -\,\textup{div}_\textup{Euclidean}
   (\lambda_{j n+1},\ldots,\lambda_{j N}),\;j=1,\ldots,q.
\end{equation}
   \item[\small{(\emph{Step 2})}] \label{page:passo 2} First we suppose that
   $\utilde\in C_c^{\infty}(W)$ is a weak solution of the equation
\begin{equation}\label{f:luxtilde_uguale_fcaptilde_piu_g_tutto_c_infinito}
   \xjttilde(\aijtilde\xitilde\utilde)=\xjttilde\fcapjtilde+g,
\end{equation}
   where $g\in C_c^{\infty}(W)$ is fixed, $\fcaptilde\in C_c^{\infty}(W,\erreq),$
   $\aijtilde\in C^{\infty}(W)\cap VMO_{\xtilde}(W);$ we prove
   that if $\textrm{supp\,} u\subset B_{\xtilde},$
   with $B_{\xtilde}$ open $C\textrm{--}\,C$ ball with radius
   $\sigma<\overline{\sigma},$ $\overline{\sigma}$ sufficiently small, 
   than hypothesys $\aijtilde\in VMO_{\xtilde}(W)$ yields
\begin{equation}\label{f:prima_stima_tutto_c_infinito_c}
   \|\xtilde \utilde\|_{L^p(B_{\xtilde})}\leq
   c\Big(\|\utilde\|_{L^p(B_{\xtilde})}+
   \|\fcaptilde\|_{L^p(B_{\xtilde},\erreq)}+
   \|\xtilde \utilde\|_{L^{p_\ast}(B_{\xtilde})}+
   \|g\|_{L^{p_\ast}(B_{\xtilde})}\Big),
\end{equation}
   where $p_\ast$ is as in Notation \eqref{notaz:p_stars},
   Section \eqref{ss:op_diff_sol_fond_par}.
   \newline
   \indent
   This step requires more care then the following ones; indeed we have to
   prove the existence of a parametrix, more precisely of an operator that,
   up to a finite sum of operators of type $0,1$ and $2,$
   behaves like a right inverse and, actually, also as a left inverse
   because of the simmetry of the matrix with entries $\aij.$
   The costruction of this parametrix, as shown in Section
   \eqref{ss:op_diff_sol_fond_par}, makes use of results in
   Section \textsc{\lq\lq Part III. Operators corresponding to free vector
   fields\rq\rq} in \cite{rothstein}, and Section
   \textsc{\lq\lq 2. Differential operators and fundamental solutions\rq\rq}
   in \cite{brmbrn2}; in particular the estimate is obtained through
   results of Section \eqref{ss:op_diff_sol_fond_par}.
   \item[\small{(\emph{Step 3})}] Now we suppose that
   $\utilde\in C^{\infty}(W),$ $\fcaptilde\in C^{\infty}(W,\erreq),$ $\aijtilde\in
   C^{\infty}(W),$ $g\in C^{\infty}(W)$ is fixed and $\utilde$ is a solution of
   \eqref{f:luxtilde_uguale_fcaptilde_piu_g_tutto_c_infinito}.
   Fix $\xi\in W,$ $0<\gamma<1$ and $0<\sigma<\overline{\sigma}.$
   Arguing as in Lemma 3.3 in \cite{brmbrn2} we choose
   $\theta\in C_c^{\infty}(W)$ such that
   $B'_{\xtilde}\prec\theta\prec B''_{\xtilde},$ where
   $B'_{\xtilde}\Subset B''_{\xtilde}\Subset W$
   are concentric balls around $\xi,$ with respective radii
   $\gamma\sigma<\sigma,$
   and such that $|\xtilde \theta|\leq\frac{c}{(1-\gamma)\sigma}.$
   Then the function
   $\theta u\in C_c^{\infty}(W)$ ia a weak solution of 
\begin{equation}\label{f:luxtilde_uguale_localizzati}
   \xjttilde(\aijtilde\xitilde(\theta \utilde))=
   \xjttilde(\aijtilde(\utilde \xitilde \theta + \theta\fcapjtilde))+
   \Big(g \theta - \aijtilde (\xitilde \utilde)(\xjtilde \theta)\Big)
\end{equation}
   and, by applying previous Step 2, the following estimate holds
\begin{equation}\label{f:stima_con_localizzata}
   \|\xtilde \utilde\|_{L^p(B'_{\xtilde})}\leq
   c\Big(\|\utilde\|_{L^p(B''_{\xtilde})}+
   \|\fcaptilde\|_{L^p(B''_{\xtilde},\erreq)}+
   \|\xtilde \utilde\|_{L^{p_\ast}(B''_{\xtilde})}+
   \|g \|_{L^{p_\ast}(B''_{\xtilde})}\Big).
\end{equation}
   \item[\small{(\emph{Step 4})}] \label{page:passo_4}
   Fixed $2<p\leq 2^{\ast},$ let us suppose that $\utilde\in W_{\xtilde}^{1,p}(W)$
   is a solution of \eqref{f:luxtilde_uguale_fcaptilde_piu_f}. Let
   $\aijtilde_h\in C^{\infty}(W)\cap VMO_{\xtilde}(W)$ $q^2$ function's sequences
   converging respectively to $\aijtilde,$ let
   $\fcaptilde_h\in C^{\infty}(W)\cap L^p(W,\erreq)$ a function's sequence
   converging in $L^p(W,\erreq)$ to $\fcaptilde,$
   and let $f_h\in C^{\infty}(W)\cap L^{p}(W)$ a function's sequence
   converging in $L^{p}(W)$ to $f.$ Let us consider the sequence of
   Dirichlet problems
\begin{equation}\label{f:probls_dirichlet}
   (D_h)\quad{
\begin{cases}
   \xjttilde(\aijtilde_h(\xitilde
   \utilde_h))=\xjttilde\fcapjtilde_h+f_h \\
   \utilde_h-\utilde\in \overset{\circ}{H}_{\xtilde}^1(W)
\end{cases}}.
\end{equation}
   For any $h=1,2,\ldots,$ let $\utilde_h\in H_{\xtilde}^1(W)$
   the unique weak solution of $(D_h).$
   Then, for any $h=1,2,\ldots,$ the function $v_h=\utilde_h-\utilde$ is
   solution of the problem
\begin{equation}\label{f:probls_dirichlet_v_con_h}
   (D_h')\quad{
\begin{cases}
   \xjttilde(\aijtilde_h(\xitilde
   v_h))=\xjttilde(\fcapjtilde_h-\fcaptilde)+(f_h-f)-
   \xjttilde((\aijtilde_h-\aijtilde)\xitilde\utilde) \\
   v_h\in \overset{\circ}{H}_{\xtilde}^1(W)
\end{cases}}.
\end{equation}
   So, for any $h=1,2,\ldots,$ we have
\begin{equation}\label{f:stima_v_con_h_lm}
   \|v_h\|_{H_{\xtilde}^1(W)}\leq
   c\Big(\|\fcaptilde_h-\fcaptilde\|_{L^2(W,\erreq)}+
   \|f_h-f\|_{L^2(W)}+\|(\aijtilde_h-\aijtilde)\xitilde\utilde\|_{L^2(W)}
   \Big)\,,
\end{equation}
   from which it follows that
   $\utilde_h\to\utilde$ in $H_{\xtilde}^1(W).$
   Let us observe now that $\utilde_h\in C^{\infty}(W).$ Indeed
   the operators $\xjtilde(\aijtilde_h\xitilde)$
   have smooth coefficients in $W,$ then can be equivalently be written
   as a sum of an operator in non divergence form and one term of less degree,
   that is of the kind $-\aijtilde_h\xjtilde\xitilde + b_h^i\xitilde,$
   with $b^i\in C^{\infty}(W).$ From the properties of convolution, and arguing as in
   Teorema 1.11 di \cite{brmbrn1}, we obtain that any $\xjtilde(\aijtilde_h\xitilde)$
   is hypoelliptic. It follows that $\utilde_h\in W_{\xtilde}^{1,p}(W)$ for any
   $1\leq p\leq\infty.$ So we can apply to each function $\utilde_h$ local
   De Giorgi--Stampacchia--Moser estimates; more precisely, when $p\geq 2,$
   because $\utilde_h\in H_{\xtilde}^{1}(W)$ is solution of the equation in
   $(D_h),$ it follows that, up to a shrinking of $W,$ there exists
   $0<\widetilde{\sigma}<\overline{\sigma},$ such that, for any $\xi\in W,$
   we can find two $C\textrm{--}\,C$ open balls, say
   $\overline{B}'_{\xtilde}$ e $\overline{B}''_{\xtilde},$
   of suitable radii $0<\ro'<\ro''<\widetilde{\sigma},$
   not depending on the given point, and a constant
   $c=c(p,\ro',\ro'',\nu,Q)$ such that
   $\|\utilde_h\|_{L^\infty(\overline{B}'_{\xtilde})}\leq c
   \big(\|\utilde_h\|_{L^2(\overline{B}''_{\xtilde})}+1\big):$
   last formula implies that $\utilde_h$ are uniformly bounded in
   $L^p(\overline{B}'_{\xtilde}).$
   Applying now \emph{Step 3}, for any point $\xi\in W$ we find two
   $C\textrm{--}\,C$ open balls, say $B'_{\xtilde}\Subset B''_{\xtilde},$
   with radii $0<\gamma\sigma<\sigma<\widetilde{\sigma}$
   such that
\begin{equation}\label{f:stima_v_con_h_lp}
   \|\xtilde \utilde_h\|_{L^p(B'_{\xtilde})}\leq
   c\Big(\|\utilde_h\|_{L^p(B''_{\xtilde})}+
   \|\fcaptilde_h\|_{L^p(B''_{\xtilde},\erreq)}+
   \|\xtilde \utilde_h\|_{L^{p_{_{\ast}}}(B''_{\xtilde})}+
   \|f_h \|_{L^{p_{_{\ast}}}(B''_{\xtilde})}\Big)
\end{equation}
   and, moreover, $u_h$ are uniformly bounded in the biggest ball.
   This fact, considering that $p_\ast<2$ and that
   $\{\utilde_h\}$ is bounded in $H_{\xtilde}^{1}(W)$ implies
   $\|\utilde_h\|_{W_{\xtilde}^{1,p}(B'_{\xtilde})}\leq$ constant for any
   $h=1,2,\ldots\,.$ So there exists a subsequence, that we keep calling
   $\{\utilde_h\},$ which weakly converges in $W_{\xtilde}^{1,p}(B'_{\xtilde})$
   to a given $\overline{u}\in W_{\xtilde}^{1,p}(B'_{\xtilde}).$
   But from
   $\utilde_h\stackrel{\tiny{W_{\xtilde}^{1,2}(B'_{\xtilde})}}
   {\rightharpoonup}\utilde$
   it follows that
   $\utilde_h\stackrel{\tiny{W_{\xtilde}^{1,p}(B'_{\xtilde})}}
   {\rightharpoonup}\utilde,$
   so $\utilde=\overline{u}$ in $B'_{\xtilde}.$ By the uniform
   boundedness principle finally it follows that, for any $2<p\leq 2^{\ast},$
\begin{align}\label{f:stima_utilde_dopo_sistema}
   \|\xtilde \utilde\|_{L^p(B'_{\xtilde})}
   &\leq
   c\Big(\|\utilde\|_{L^p(B''_{\xtilde})}+
   \|\fcaptilde\|_{L^p(B''_{\xtilde},\erreq)}+
   \|\xtilde \utilde\|_{L^{p_\ast}(B''_{\xtilde})}+
   \|f \|_{L^{p_\ast}(B''_{\xtilde})}\Big)
   \notag\\
   &\leq
   c\Big(\|\utilde\|_{L^p(B''_{\xtilde})}+
   \|\fcaptilde\|_{L^p(B''_{\xtilde},\erreq)}+
   \|\xtilde\utilde\|_{L^{p_\ast}(B''_{\xtilde})}\Big)
   \notag\\
   &\leq
   c\Big(\|\utilde\|_{L^p(B''_{\xtilde})}+
   \|\fcaptilde\|_{L^p(B''_{\xtilde},\erreq)}+
   \|\xtilde\utilde\|_{L^{2}(B''_{\xtilde})}\Big).
   \notag\\
\end{align}
   \item[\small{(\emph{Step 5})}] Here we employ the recursive tecnique as in
   \cite{difazio}. More precisely, by applying repetitively \emph{Step 3} and
   \emph{Step 4}, we will show the existence of $L^p$ estimate when the number $Q$
   \lq\lq grows up\rq\rq with $p.$ More precisely we have
   $]2,\infty[=]2,2^{\ast}]\cup]2^{\ast},2^{\ast\ast}]\cup\ldots\,$
   where, for $2\leq q < Q,$ $q^{\ast}$ is defined as in Notation \eqref{notaz:p_stars}. 
   In particolar, for $q=2,$ $q^{\ast\ast}$ is defined only when $Q>4,$ and so on.
   Let us suppose that $\utilde\in W_{\xtilde}^{1,p}(W)$ is a solution of
   \eqref{f:luxtilde_uguale_fcaptilde_piu_f}. Then \emph{Step 4} implies that
   if $2<p\leq 2^{\ast}$ and $0<\sigma<\overline{\sigma}$
   is sufficiently small, for any $\xi\in W,$ and relatively to balls
   $B'_{\xtilde}\Subset B''_{\xtilde}$ around $\xi\in W$ with radii
   $0<\gamma\sigma<\sigma,$ \eqref{f:stima_utilde_dopo_sistema} holds.
   Fixed now $2^\ast<p\leq 2^{\ast\ast}$ and repeat again arguments in
   \emph{Step 3} and \emph{Step 4}. Choose a function
   $\theta\in C_c^{\infty}(B'_{\xtilde}),$ where $B'_{\xtilde}$ is,
   for any fixed $\xi,$ the open $C\textrm{--}\,C$
   ball around $\xi$ with radius $\gamma\sigma;$
   more precisely we choose this ball such that
   $\overline{B}_{\xtilde}\prec\theta\prec B'_{\xtilde},$ where
   $\overline{B}_{\xtilde}\Subset B'_{\xtilde}$
   are two open $C\textrm{--}\,C$ concentric balls around $\xi,$ with radii
   $0<\gamma^2\sigma<\gamma\sigma,$ and such that
   $|\xtilde \theta|\leq\frac{c}{(1-\gamma)\gamma\sigma}.$
   Then
\begin{align}\label{f:prima_stima_utilde_passo_5}
   \|\xtilde\utilde\|_{L^p(\overline{B}_{\xtilde})}
   &\leq
   c\Big(\|\utilde\|_{L^p(B'_{\xtilde})}+
   \|\fcaptilde\|_{L^p(B'_{\xtilde},\erreq)}+
   \|\xtilde\utilde\|_{L^{p_\ast}(B'_{\xtilde})}\Big)
   \notag\\
   &\leq
   c\Big(\|\utilde\|_{L^p(B'_{\xtilde})}+
   \|\fcaptilde\|_{L^p(B'_{\xtilde},\erreq)}+
   \|\xtilde\utilde\|_{L^{2^\ast}(B'_{\xtilde})}\Big),
   \notag\\
   &\leq
   c\Big(\|\utilde\|_{L^p(B''_{\xtilde})}+
   \|\fcaptilde\|_{L^p(B''_{\xtilde},\erreq)}+
   \|\utilde\|_{L^{2^\ast}(B''_{\xtilde})}+
   \|\fcaptilde\|_{L^{2^\ast}(B''_{\xtilde},\erreq)}+
   \|\xtilde\utilde\|_{L^{2}(B''_{\xtilde})}\Big)
   \notag\\
   &\leq
   c\Big(\|\utilde\|_{L^p(B''_{\xtilde})}+
   \|\fcaptilde\|_{L^p(B''_{\xtilde},\erreq)}+
   \|\xtilde\utilde\|_{L^{2}(B''_{\xtilde})}\Big).
\end{align}
   So, if $p>2$ let $h\in\mathbb{N}$ be such that
   $2^{\overbrace{\ast\ast\ldots\ast}^{h-1}}<p\leq
   2^{\overbrace{\ast\ast\ldots\ast}^{h}}.$ If
   $\gamma=(\frac{1}{2})^{\frac{1}{h}},$ by reitering the argument we obtain
\begin{equation}\label{f:seconda_stima_utilde_passo_5}
   \|\xtilde\utilde\|_{L^p(B'_{\xtilde})}\leq
   c\Big(\|\utilde\|_{L^p(B''_{\xtilde})}+
   \|\fcaptilde\|_{L^p(B''_{\xtilde},\erreq)}+
   \|\xtilde\utilde\|_{L^{2}(B''_{\xtilde})}\Big)
\end{equation}
   where $B'_{\xtilde}\Subset B''_{\xtilde}$ are two open
   $C\textrm{--}\,C$ balls of radii $0<\frac{\sigma}{2}<\sigma,$
   for any $0<\sigma<\widetilde{\sigma}.$
   \item[\small{(\emph{Step 6})}] So far we have obtained that
   if $2\leq p <\infty$ and $\utilde\in W_{\xtilde}^{1,p}(W)$ is a solution
   of \eqref{f:luxtilde_uguale_fcaptilde_piu_f}, then the estimate
   \eqref{f:seconda_stima_utilde_passo_5} holds. Observe now that
   in particular $\utilde\in W_{\xtilde}^{1,2}(W)$ is a solution of
   \eqref{f:luxtilde_uguale_fcaptilde_piu_f}: then classical
   $L^2$ theory yields the existence of $0<\ro<\widetilde{\sigma}$
   such that for any $2\sigma<\ro,$
   if $B'''_{\xtilde}$ is the open ball of radius $2\sigma$
   concentric with $B''_{\xtilde},$ it results
\begin{align*}
   \|\xtilde\utilde\|_{L^{2}(B''_{\xtilde})}&\leq
   c\Big(\|\utilde\|_{L^2(B'''_{\xtilde})}+
   \|\fcaptilde\|_{L^2(B'''_{\xtilde},\erreq)}\Big)\\
   &\leq c\Big(\|\utilde\|_{L^p(B'''_{\xtilde})}+
   \|\fcaptilde\|_{L^p(B'''_{\xtilde},\erreq)}\Big)\\
\end{align*}
   because $2<p.$
   This last estimate, jointly with \eqref{f:seconda_stima_utilde_passo_5}
   implies \eqref{f:stima_finale_palla_palla_doppia}.
   \item[\small{(\emph{Step 7})}] Noting now that $\sigma$ is as
   small as we need, in
   \eqref{f:stima_finale_palla_palla_doppia}, taking into account
   the local equivalence of $d_{\widetilde{X}}$ and the Euclidean metric
   we can assume that \eqref{f:stima_finale_palla_palla_doppia} holds
   with two Euclidean balls; then the estimate clearly holds between
   due open sets of the kind $B\times I,B'\times I,$ with
   $B\Subset B'\subseteq\Omega$ open
   Euclidean balls small enough and $I$ open rectangle in $\errek,$ previously
   fixed small enogh. Then through Proposition 1.4 of \cite{frsese} (applied with
   the weight function $w\equiv 1$) and Remark
   \eqref{o:come_riportare_la_stima_indietro}, we obtain the existence
   of an absolute costant such that, for any
   $u\in W_{X}^{1,p}(W)$ solution of \eqref{f:lux_uguale_fcap},
   the following estimate holds
\begin{align}\label{f:formula_pre_finale}
   \|u\|_{W_{X}^{1,p}(B)}&
   \leq c\Big(\|u\|_{L^p(B')}+\|F\|_{L^p(B',\erreq)}\Big),
\end{align}
   and so the thesis.
\end{itemize}
   \noindent
   \indent
   \vskip10pt
   Finally let us sketch some of the proofs.
   \newline
   \noindent
   \indent
\begin{itemize}
   \item[\small{(\emph{Step 1 - Some Remarks})}] If $u\in S_X^{1,p}(\Omega),$
   then $\utilde\in W_{\xtilde}^{1,p}(W).$ Indeed, by Proposition 1.4
   in \cite{frsese}, we can assume that $\utilde\in C^\infty(W).$
   Then $\mathcal{D}(B)\otimes\mathcal{D}(I)$ is dense into
   $\mathcal{D}(W)$ and, in particular, we can consider
   test functions of the kind
   $\Phi(\xi)\equiv\phi(x)\psi(t)$ for any $\xi=(x,t)\in B\times I,$
   with arbitrary $\phi$ e $\psi$ in $\mathcal{D}(B)$ and
   $\mathcal{D}(I)$. So, in
   \eqref{f:lux_uguale_fcap_per_esteso}, it suffices to multiply by
   $\psi(t)$ and integrate over $I,$ in order to apply Remark
   \eqref{o:come_riportare_la_stima_indietro}. In particular
   $\widetilde{X_j}(u\circ \pi)=X_ju\circ \pi$ holds for any
   $u\in S_X^{1,p}(\Omega).$ Analogous arguments hold for
   \eqref{f:luxtilde_uguale_fcaptilde_piu_f}.
   \item[\small{(\emph{Step 2 - Proof})}]
   We refer to notations and results in Section \eqref{ss:op_diff_sol_fond_par}.
   Fix $\xi_0\in W$ and let us consider the operator
   $\widetilde{L}_0\equiv$
   $\widetilde{L}_{\xi_0}=\xjttilde(\aijtilde(\xi_o)\xitilde).$
   We have to construct two partial inverse $D_0$ and $S_0$ of
   $\widetilde{L}_0$ of type 2, frozen at $\xi_0.$ We have
   $\widetilde{L}=\widetilde{L}_0+(\widetilde{L}-\widetilde{L}_0);$
   then we can apply to both member of the equation a suitable operator
   of type $1$ and, finally, after some calculations, we can free $\xi_0$
   having so a representation of $\xitilde u,$ so to apply Theorems
   \eqref{t:limitatez_int_sing} and \eqref{cor:i_due_tipi_vanno_a_star_basso}
   to obtain \eqref{f:prima_stima_tutto_c_infinito_c}.
   Precisely, let us fix $a\in C_c^\infty(W)$ and $\xi_0\in W.$ We are going to prove 
   that there exist two frozen operators at $\xi_0,$ $D_{0,2}$ and $S_{0,2}$ of
   type 2, and a finite number of operators frozen at $\xi_0$
   (depending only on $\lux\,$)
   $D_{0,1}^{ij,\,h}$ e $D_{0,2}^{ij,\,k},$ $h\in I,\;k\in J$ with
   $|I|,|J|\leq\textup{absolute constant},$
   which are Riesz potential of type $1$ and $2,$ such that,
   for any test $\utilde$ in $W$ it results
\begin{equation}\label{f:parametrice_destra}
   \widetilde{L}_0D_{0,2}\utilde=-au+
   \sum_{h\in I}\sum_{i,j=1}^{q}\aijtilde(\xi_0)D_{0,1}^{ij,\,h}\utilde+
   \sum_{k\in
   K}\sum_{i,j=1}^{q}\aijtilde(\xi_0)D_{0,2}^{ij,\,k}\utilde,
\end{equation}
\begin{equation}\label{f:parametrice_sinistra}
   S_{0,2}\widetilde{L}_0\utilde=-au+
   \sum_{h\in I}\sum_{i,j=1}^{q}\aijtilde(\xi_0)S_{0,1}^{ij,\,h}\utilde+
   \sum_{k\in K}\sum_{i,j=1}^{q}\aijtilde(\xi_0)S_{0,2}^{ij,\,k}\utilde.
\end{equation}
   To this aim, let us consider the fundamental solution $\Gamma_0$ ensured by
   Folland's theorem, of the invariant operator
   $\ajitilde(\xi_0)Y_jY_i$ associated to the non divergence form operator
   $\ajitilde(\xi_0)\xjtilde\xitilde$ (see Theorem \eqref{t:teor_folland_sol_fond}).
   Let us consider now the operator that in \cite{brmbrn2} is adapted
   to the one of \cite{rothstein}; more precisely, fixed a test
   $b$ in $W$ such that $\textup{supp }a\subset \{b=1\},$
   according to notations of Theorem \eqref{t:teta}, we set,
   for any $\utilde$ test in $W,$ and for any $\xi\in W,$
\begin{equation}\label{parametrice}
   D_{0,2}\utilde(\xi)=\frac{a(\xi)}{\omega(\xi)}\int_W
   \Gamma(\tetaex)b(\eta)\utilde(\eta)d\eta.
\end{equation}
   Applying Theorem \eqref{t:appr}, for any $i=1,\ldots,q$ we have
\begin{align}\label{f:xconi_applicato_a_d(0,2)}
   \xitilde D_{0,2}\utilde(\xi)&=
   \xitilde
   \left(
   \frac{a(\xi)}{\omega(\xi)}
   \right)
   \int_W\Gamma(\tetaex)b(\eta)\utilde(\eta)d\eta\notag\\
   &+\frac{a(\xi)}{\omega(\xi)}
   \int_W
   \left[
   Y_i\Gamma_0(\tetaex)+R_i^\xi\Gamma_0(\tetaex)
   \right]
   b(\eta)\utilde(\eta)d\eta.
\end{align}
   Now, according to Theorem \eqref{t:liber},
   if $c_j=(a_{j\,1},\ldots,a_{j\,n},\lambda_{j\,n+1},\ldots,\lambda_{j\,N})$
   denote entries of the vector field $\xjtilde,$ it is
   $\xjttilde=-\xjtilde+m_j$ where $m_j(\xi)=$
   $(\textup{div}_{\textup{euclidea}}\,c_j)(\xi),$ for any $\xi\in W.$ 
   So, arguing as in \cite{brmbrn2}, we have
\begin{align*}
   \widetilde{L}_0D_{0,2}\utilde(\xi)=
   -\xjtilde
   &\bigg[
   \ajitilde(\xi_0)\xitilde
   \left(
   \frac{a(\xi)}{\omega(\xi)}
   \right)
   \int_W\Gamma_0(\tetaex)b(\eta)\utilde(\eta)d\eta
   \notag\\
   &+\ajitilde(\xi_0)\frac{a(\xi)}{\omega(\xi)}
   \int_W
   \left[
   Y_i\Gamma_0(\tetaex)+R_i^\xi\Gamma_0(\tetaex)
   \right]
   b(\eta)\utilde(\eta)d\eta
   \bigg]\\
   +m_j&(\xi)
   \bigg[
   \ajitilde(\xi_0)\xitilde
   \left(
   \frac{a(\xi)}{\omega(\xi)}
   \right)
   \int_W\Gamma_0(\tetaex)b(\eta)\utilde(\eta)d\eta\\
   &+\ajitilde(\xi_0)\frac{a(\xi)}{\omega(\xi)}
   \int_W
   \left[
   Y_i\Gamma_0(\tetaex)+R_i^\xi\Gamma_0(\tetaex)
   \right]b(\eta)\utilde(\eta)d\eta
   \bigg]=
\end{align*}
\begin{align}\label{f:ltilde0_applicato_a_d(0,0)}
   =\ajitilde(\xi_0)
   &\bigg[
   \xjtilde\xitilde
   \left(
   \frac{a(\xi)}{\omega(\xi)}
   \right)
   \int_W\Gamma_0(\tetaex)b(\eta)\utilde(\eta)d\eta\notag\\
   &+\xitilde
   \left(
   \frac{a(\xi)}{\omega(\xi)}
   \right)
   \int_W
   \left[
   Y_j\Gamma_0(\tetaex)+R_j^\xi\Gamma_0(\tetaex)
   \right]
   b(\eta)\utilde(\eta)d\eta\notag\\
   &+\xjtilde
   \left(
   \frac{a(\xi)}{\omega(\xi)}
   \right)
   \int_W
   \left[
   Y_i\Gamma_0(\tetaex)+R_i^\xi\Gamma_0(\tetaex)
   \right]
   b(\eta)\utilde(\eta)d\eta\notag\\
   &+\frac{a(\xi)}{\omega(\xi)}
   P.V.\int_W
   \Big[
   Y_jY_i\Gamma_0(\tetaex)+Y_jR_i^\xi\Gamma_0(\tetaex)+\notag\\
   &\qquad\qquad\qquad\qquad\qquad\qquad
   +R_j^\xi Y_i\Gamma_0(\tetaex)+R_j^\xi R_i^\xi\Gamma_0(\tetaex)
   \Big]
   b(\eta)\utilde(\eta)d\eta
   \bigg]\notag\\
   +m_j(\xi)&
   \bigg[
   \ajitilde(\xi_0)\xitilde
   \left(
   \frac{a(\xi)}{\omega(\xi)}
   \right)
   \int_W\Gamma_0(\tetaex)b(\eta)\utilde(\eta)d\eta\notag\\
   &+\ajitilde(\xi_0)\frac{a(\xi)}{\omega(\xi)}
   \int_W
   \left[
   Y_i\Gamma_0(\tetaex)+R_i^\xi\Gamma_0(\tetaex)
   \right]b(\eta)\utilde(\eta)d\eta
   \bigg]\notag\\
   =-au+
   &\sum_{h\in I}\sum_{i,j=1}^{q}\aijtilde(\xi_0)D_{0,1}^{ij,\,h}\utilde+
   \sum_{k\in
   K}\sum_{i,j=1}^{q}\aijtilde(\xi_0)D_{0,2}^{ij,\,k}\utilde.
\end{align}
   \indent
   Finally, by transposition of the matrix $(\aijtilde)_{i,j=1,\ldots,q}$
   we obtain the desired formula.
   In particular, fixed any test $\utilde,$ for any test function
   $a$ such that
   $\textup{supp }\utilde\subset \{a=-1\},$ we can write
   $\xtilde_m\utilde.$ More precisely, by applying
   $\xtilde_m$ to both member of
   \eqref{f:parametrice_sinistra},
   by Theorem \eqref{t:xtilde_con_operatore_scambiano}, for any $m=1,\ldots,q$
   it is
\begin{align}\label{f:xcon_m_para_sinistra}
   \xtilde_m S_{0,2} \widetilde{L}_0 u=\xtilde_m u
   &+\sum_{h\in I}\sum_{i,j=1}^{q}\aijtilde(\xi_0)
   \left[\sum_{l=1}^q T_{0,1}^{ij,\,hl}\widetilde{X_l}
   +T^{ij,\,h}_{0,1}\right]\utilde\notag\\
   &+\sum_{k\in K}\sum_{i,j=1}^{q}\aijtilde(\xi_0)
   \left[\sum_{l=1}^q T_{0,2}^{ij,\,kl}\widetilde{X_l}
   +T^{ij,\,k}_{0,2}\right]\utilde.
\end{align}
   \newline
   \indent
   It follows that, thanks to Lemma \eqref{l:tipo_i_tipo_i-1},
   $T_{0,1}\equiv\xtilde_m S_{0,2},$ is an operator of type 1 frozen
   at $\xi_0.$ So, let $\utilde$ be a solution of
   \eqref{f:luxtilde_uguale_fcaptilde_piu_g_tutto_c_infinito}.
   By applying $T_{0,1}$ to
\begin{align}\label{f:aggiungi_sottrai}
   \widetilde{L}_0\utilde&=\widetilde{L}\utilde+
   (\widetilde{L}_0-\widetilde{L})\utilde\notag\\
   &=\xjttilde\fcapjtilde+g+
   \xjttilde((\aijtilde(\xi_0)-\aijtilde)\xitilde\utilde)
\end{align}
   it is
\begin{align}\label{f:t(0,1)_applicato_alla_aggiungi_sottrai}
   T_{0,1}\widetilde{L}_0\utilde&
   =T_{0,1}\xjttilde\fcapjtilde+T_{0,1}g+
   T_{0,1}
   \left[
   \xjttilde((\aijtilde(\xi_0)-\aijtilde)\xitilde\utilde)
   \right].
\end{align}
   Then, putting \eqref{f:t(0,1)_applicato_alla_aggiungi_sottrai}
   into \eqref{f:xcon_m_para_sinistra} and getting $\xtilde_m\utilde,$
   we have
\begin{align}\label{f:xitilde_con_m}
   \xtilde_m\utilde&=T_{0,1}\xjttilde\fcapjtilde+T_{0,1}g\notag\\
   &+T_{0,1}
   \left[
   \xjttilde((\aijtilde(\xi_0)-\aijtilde)\xitilde\utilde)
   \right]\notag\\
   &-\sum_{h\in I}\sum_{i,j=1}^{q}\aijtilde(\xi_0)
   \left[
   \sum_{l=1}^q T_{0,1}^{ij,\,hl}\widetilde{X_l}
   +T^{ij,\,l}_{0,1}
   \right]
   \utilde\notag\\
   &-\sum_{k\in K}\sum_{i,j=1}^{q}\aijtilde(\xi_0)
   \left[
   \sum_{l=1}^q T_{0,2}^{ij,\,kl}\widetilde{X_l}
   +T^{ij,\,k}_{0,2}
   \right]\utilde.
\end{align}
   Finally, through the definition of transposed operator and
   Lemma \eqref{l:tipo_i_tipo_i-1}, it follows that for $q$
   suitable operators $T_{0,0}^j$ frozen at $\xi_0$ of type $0,$ 
   \eqref{f:xitilde_con_m} becomes
\begin{align}\label{f:xitilde_con_m_finale}
   \xitilde_m\utilde&=\sum_{j=1}^q T_{0,0}^j\fcapjtilde+T_{0,1}g\notag\\
   &+\sum_{j=1}^q T_{0,0}^j
   ((\aijtilde(\xi_0)-\aijtilde)\xitilde\utilde)\notag\\
   &-\sum_{h\in I}\sum_{i,j=1}^{q}\aijtilde(\xi_0)
   \left[
   \sum_{l=1}^q T_{0,1}^{ij,\,hl}\widetilde{X_l}
   +T^{ij,\,l}_{0,1}
   \right]
   \utilde\notag\\
   &-\sum_{k\in K}\sum_{i,j=1}^{q}\aijtilde(\xi_0)
   \left[
   \sum_{l=1}^q T_{0,2}^{ij,\,kl}\widetilde{X_l}
   +T^{ij,\,k}_{0,2}
   \right]\utilde.
\end{align}
   \indent
   From this representation and the arbitrariness of $\xi_0\in W,$ 
   \eqref{f:prima_stima_tutto_c_infinito_c} follows immediately by
   applyng Theorem \eqref{t:limitatez_int_sing}
   and Corollary \eqref{cor:i_due_tipi_vanno_a_star_basso}.
\end{itemize}
   \vskip15pt
   \indent
   This concludes the sketch of the proof of Theorem \eqref{teor:regolarita}.$\B$

\end{document}